\documentclass[11pt, a4paper]{amsart}
\usepackage{verbatim}
\usepackage[dvips]{color}
\usepackage{amssymb}
\usepackage[active]{srcltx}
\usepackage{latexsym}
\usepackage{amsmath}
 \usepackage{float}
\usepackage{mathrsfs} 
\usepackage[draft]{fixme}
\usepackage{graphicx, epstopdf, float}

\allowdisplaybreaks

\theoremstyle{definition}

\floatstyle{ruled} \restylefloat{figure} \restylefloat{table}

\numberwithin{equation}{section}

\newcommand{\beq}{\begin{equation}}
\newcommand{\bea}[1]{\begin{array}{#1} }
\newcommand{\eeq}{ \end{equation}}
\newcommand{\ea}{ \end{array}}

\def\mean#1{\mathchoice%
          {\mathop{\kern 0.2em\vrule width 0.6em height 0.69678ex depth -0.58065ex
                  \kern -0.8em \intop}\nolimits_{\kern -0.4em#1}}%
          {\mathop{\kern 0.1em\vrule width 0.5em height 0.69678ex depth -0.60387ex
                  \kern -0.6em \intop}\nolimits_{#1}}%
          {\mathop{\kern 0.1em\vrule width 0.5em height 0.69678ex
              depth -0.60387ex
                  \kern -0.6em \intop}\nolimits_{#1}}%
          {\mathop{\kern 0.1em\vrule width 0.5em height 0.69678ex depth -0.60387ex
                  \kern -0.6em \intop}\nolimits_{#1}}}

\def\vintslides_#1{\mathchoice%
          {\mathop{\kern 0.1em\vrule width 0.5em height 0.697ex depth -0.581ex
                  \kern -0.6em \intop}\nolimits_{\kern -0.4em#1}}%
          {\mathop{\kern 0.1em\vrule width 0.3em height 0.697ex depth -0.604ex
                  \kern -0.4em \intop}\nolimits_{#1}}%
          {\mathop{\kern 0.1em\vrule width 0.3em height 0.697ex depth -0.604ex
                  \kern -0.4em \intop}\nolimits_{#1}}%
          {\mathop{\kern 0.1em\vrule width 0.3em height 0.697ex depth -0.604ex
                  \kern -0.4em \intop}\nolimits_{#1}}}

\newcommand{\aveint}[2]{\mathchoice%
          {\mathop{\kern 0.2em\vrule width 0.6em height 0.69678ex depth -0.58065ex
                  \kern -0.8em \intop}\nolimits_{\kern -0.45em#1}^{#2}}%
          {\mathop{\kern 0.1em\vrule width 0.5em height 0.69678ex depth -0.60387ex
                  \kern -0.6em \intop}\nolimits_{#1}^{#2}}%
          {\mathop{\kern 0.1em\vrule width 0.5em height 0.69678ex depth -0.60387ex
                  \kern -0.6em \intop}\nolimits_{#1}^{#2}}%
          {\mathop{\kern 0.1em\vrule width 0.5em height 0.69678ex depth -0.60387ex
                  \kern -0.6em \intop}\nolimits_{#1}^{#2}}}

\def\eqn#1$$#2$${\begin{equation} \label#1#2\end{equation}}
\def\charfn_#1{{\raise1.2pt\hbox{$\chi
_{\kern-1pt\lower3pt\hbox{{$\scriptstyle#1$}}}$}}}

\def\qq1{q_*}
\def\q2{q_{**}}

\newdimen\vintbar
\def\vint{-\kern-\vintbar\int}

\def\0{\boldsymbol 0}

\newtoks\by
\newtoks\paper
\newtoks\book
\newtoks\jour
\newtoks\yr
\newtoks\pages
\newtoks\vol
\newtoks\publ

\def\name[#1, #2]{#1 #2}
\def\ota{{\hbox{\bf ???}}}
\def\cLear{\by = \ota\paper = \ota\book = \ota\jour = \ota\yr = \ota
\pages = \ota\vol = \ota\publ = \ota}
\def\endpaper{\the\by, \textit{\the\paper},
{\the\jour} \textbf{\the\vol} (\the\yr), \the\pages.\cLear}
\def\endbook{\the\by, \textit{\the\book},
\the\publ, \the\yr.\cLear}
\def\endpap{\the\by, \textit{\the\paper}, \the\jour.\cLear}
\def\endproc{\the\by, \textit{\the\paper}, \the\book, \the\publ,
\the\yr, \the\pages.\cLear}

\begin{document}


\title{A framework for the Modeling of Order Book\\ Dynamics based on Event Sizes}\thanks{This research has been supported by a grant from Riksbankens Jubileumsfond (P10-0113:1).} \thanks{The empirical analysis in this paper is carried out using high frequency data supplied by NasdaqOMX. The first author like to in particularly thank Petter Dahlstr{\"o}m and Bj{\"o}rn Hertzberg, both at  NasdaqOMX in Stockholm, for supplying the data and for valuable discussions concerning high frequency trading.}

 \address{Kaj Nystr\"{o}m\\Department of Mathematics, Uppsala University\\
S-751 06 Uppsala, Sweden}
\email{kaj.nystrom@math.uu.se}

\address{Sidi Mohamed Ould Aly\\Department of Mathematics, Uppsala University\\
S-751 06 Uppsala, Sweden}
\email{souldaly@math.uu.se}


\author{Kaj Nystr{\"o}m and Sidi Mohamed Ould Aly}

\maketitle

\noindent
\begin{abstract}
\noindent  We propose a modeling framework for the dynamics of a reduced form order book in event time and based on event sizes. Our framework for the order book is influenced by \cite{ContLarrard12}, but compared to \cite{ContLarrard12} we allow the best bid ask spread to be larger than one tick. Based
on the modeling assumption that the best bid as well as the best ask price can only move by at most one tick (up or down), when an event occurs, we show that
the dynamics of this simplified order book is completely described by a non-linear transformation of  two processes $(X,Y)$. A key challenge in the modeling is
the empirical fact that the high frequency order flow is strongly autocorrelated, a fact we have to deal with in the modeling of $(X,Y)$. The core of our framework is a semi linear regression type model for $(X,Y)$, influence by  more classical ARMA and AR models, and one key degree of freedom is the
potentially non-linear basis functions used in  the regression. We use results from the theory of random iterative function systems to understand issues concerning stationarity and ergodicity in our models. We show how to rapidly calibrate the model by inverting block Toeplitz matrices in an efficient way. All components are worked through and explain in an application and the predictability of the model for order flows and price moves are analyzed in the context of a high frequency dataset.

\medskip

\noindent
2000  {\em Mathematics Subject Classification.}
\noindent

\medskip

\end{abstract}

\section{Introduction}
\noindent
 In recent years, numerous electronic trading systems have been designed to replace or complement traditional markets. In electronic order-driven markets, participants may submit limit order (or cancel an existing limit order), by specifying whether they wish to buy or sell, the amount (volume) and the desired price. Limit orders wait in a queue to be canceled or executed (when a sell order is matched against one or more buy limit orders). All outstanding limit orders are aggregated in a limit order book which is available for market participant. The order book at a given instant of time is the list of all buy and sell limit orders (that have not been executed/canceled) with their corresponding price and volume. The price dynamics is therefore the result of the interplay between anonymous traders who place orders in the limit order book. Nowadays, more than half of the  stock exchanges in the world are order driven and the submissions of limit orders  have massively increased. The latter makes it extremely hard to build a model for the order book and its dynamics which is realistic and still amenable to rigorous quantitative analysis.

There is now a growing literature on limits order books and attempts to modeling their dynamics. Several authors have proposed models based on self-exciting point processes (see e.g. \cite{Sahalia11}, \cite{Bauwens09}) focusing on the modeling of the frequency at which the different types of orders arrive. In other works (e.g. \cite{Parlour98}, \cite{Foucault99}) equilibrium models for limit order markets are proposed, while another category of models is based on some dynamic expected utility maximization models (e.g. \cite{Parlour98}, \cite{Rous09}). More recently,  Cont et al \cite{ContStoikov10} considered the limit order book as a queuing system where the frequencies at which orders are submitted  and canceled, all orders are assumed to be of unit size,  are modeled as independent Poisson processes. Although all these models are theoretically attractive they have some difficulties describing the reality of the order book dynamics and main drawbacks of the models include Markov type assumption concerning the order flow for some and the difficulty of  calibrating the parameters of the model for others. However the main reason that these models produce order book dynamics which deviates from what is observed in real data
stems from the fact that it is difficult, or even impossible, to predict the arrival times of orders due to the inhomogeneity of these arrival times. This is a problem even for self-exciting Hawkes processes \cite{hawkes71} where the main difficulty is to accurately and properly calibrate the kernel defining
the range of influence of previous events on the current intensity of events.

Based on the main drawback of the above models an idea is to model order book dynamics using a clock different from chronological time and in this paper we propose a modeling framework for the order book and its  dynamics in \textit{event-based time}, i.e. we are only interested in modeling  the order book as it is changed by the occurrence  of an event and by definition an event is the submission of an order, market, limit or cancelation, which changes the order book. In particular, in the proposed framework we make no reference to chronological time
and the specific times at which the events occur. Instead, for us event time is simply an index,  set to 0 at midnight, which represents an enumeration of the events and which increases by one every time a new order is submitted to the order book. The idea to use a different clock has already appeared in the literature, e.g. see Mandelbrot and Taylor \cite{Mandelbrot67},  Clark  \cite{Clark73},  and  more recently in An{\'e} and Geman \cite{ane2000} where it is showed that the normality of asset return can be recovered through a stochastic time change.

 Our framework for the order book is influenced by Cont and De Larrard  \cite{ContLarrard12} and as in \cite{ContLarrard12} we represent
 a reduced form of the state of the order book, right after event $j$, by $(S^b_j,S^a_j,Q^b_j,Q^a_j)$ where
 \begin{eqnarray}\label{not1}
  &&\mbox{the best bid price $S^b_j$ right after event $j$},\notag\\
   &&\mbox{the best ask price $S^a_j$ right after event $j$},\notag\\
    &&\mbox{the size of the best bid queue $Q^b_j$ right after event $j$},\notag\\
     &&\mbox{the size of the best ask queue $Q^a_j$ right after event $j$.}
      \end{eqnarray}
      The state of this reduced form of the order book changes when an order (market order, limit order or cancelation) effects the queues at the  best bid/ask. When the best bid queue (or best ask queue) is depleted the price moves one (or more) step(s) down (or up) and the best bid ask spread widens immediately to more than one tick. In \cite{ContLarrard12} the authors assume that after such an event the flow of limit buy/sell orders quickly fills the gap and that the spread instantly reduces to one tick again. In particular, we emphasize that in \cite{ContLarrard12} it is assumed that the second step occurs instantly and hence that the spread, i.e.   distance between the best bid and the best ask, always is one tick. Based on this assumption the process $Q_j = (Q^b_j, Q^a_j)$ can be constructed from the order flow process, $O^F$, until the next price move, as
\begin{equation}\label{order_flow}
Q_j =(Q^b_j, Q^a_j)= O^F_{j_0, j} :=  Q_{j_0} + \left( \sum_{i= j_0+ 1}^{j } V^b_i, \sum_{i= j_0+1}^{j}V^a_i \right), ~~Êj>j_0.
\end{equation}
Here
 \begin{eqnarray}\label{not2}
  &&\mbox{$V^b_i$ is the change of the size of the bid queue occurring at event number $i$},\notag\\\
   &&\mbox{$V^a_i$ is the change of the size of the ask queue occurring at event number $i$}.
      \end{eqnarray}
Note that $V^b_i$($V^a_i$) is positive if event $i$ is a limit order and negative if it is a market order or cancelation and that $V^b_i$ and $V^a_i$ can not simultaneously
 be different from zero as one event can not change the queues on both sides of the order book. Whenever a component of $O^F_{j_0, .}$ hits zero there is a price change: the price moves one step (tick) up or one step down, the spread widens to more than one tick and $Q$ deviates temporarily from the order flow process (as \textit{$Q^b$ and $Q^a$ are both strictly positive, by definition of the best bid and best ask price}). In \cite{ContLarrard12} it is assumed, when a component of $O^F_{j_0, .}$ hits zero, that the order book immediately jumps to the next step where the spread is again one tick and in \cite{ContLarrard12} the queue sizes are updated according to some distribution functions.

In this paper we consider the reduced form of the order book briefly detailed above,  but compared to \cite{ContLarrard12}
 \begin{eqnarray}\label{a1}\mbox{we allow the best bid ask spread to widen to more than one tick.}
  \end{eqnarray}
Our main modeling assumption is that
   \begin{eqnarray}\label{a2}\mbox{the best bid and ask prices can only change by one tick at each event}.
  \end{eqnarray}
We emphasize that this is a very week assumption which may even be considered as a fact in liquid markets. In particular, for our data this assumption is satisfied in 100\% of all cases (see Remark~\ref{model_assumption} below). Based on the assumption in \eqref{a2} we show that the dynamics of this simplified order book can be fully described by only two processes $(X,Y)$. The first of these processes, $X$, is the (signed) size of the new order and the second processes $Y$  is the size of the queue on the side of the order book which is changed by $X$. $X$ and $Y$ are a priori dependent and autocorrelated and the challenge is then to model these processes in a consistent way. A contribution of the paper
is that apart from the assumption in \eqref{a2}, we make
   \begin{eqnarray}\label{a3}\mbox{no additional assumptions concerning the order book dynamics.}
  \end{eqnarray}
  Instead we propose a modeling framework for the order book dynamics, through $(X,Y)$, which is completely data driven and contains no adhoc assumptions.

The rest of the paper is organized as follows. In section 2 we show how the modeling of $(Q^b_j,Q^a_j,V^b_j,V^a_j)$ is reduced to the modeling of
 processes  $(X, Y)$. We also show how $(S^b,S^a)$ is derived given $(Q^b,Q^a,V^b,V^a)$. In section 3, which is of empirical nature, we empirically describe
 the spread process found in high frequency data, the time to queue depletion as well as the statistical properties of the empirically observed
 $(X, Y)$. Section 3 serves as an empirical basis for our modeling framework proposed and outlined in section 4 and section 5. While we in section 4 give a short background on ARMA and AR models we in section 5 outline a semi linear model which is the core of our framework. An interesting aspect here is that
 we use results from the theory of random iterative function systems to understand issues concerning stationarity and ergodicity in our models. Section 6 is devoted to an application of our framework
 to the predication of the order flow and price moves in a the context of high frequency data. In section 7 we give a brief discussion and/of future research. Finally, in an appendix the inversion of Toeplitz matrices is discussed.

\section{Order Flow and Price Dynamics}\label{sec2}
\noindent In the following we outline the reduced form of the order book considered in this paper. Recall the notation $(S^b_j,S^a_j,Q^b_j,Q^a_j,V^b_j,V^a_j)$
introduced in the introduction. In the following we use the notation that $1_{  \textrm{event} }=1$ if event is true and $1_{  \textrm{event} }=0$ if event is not true.  Note that the processes $V^{a}_j$ and $ V^{b}_j$ are given in terms of $S^{b},~ S^{a},~Q^{b},~Q^{a},$ through
 \begin{eqnarray}\label{v_0}
 V^{b}_j& = &\left( Q^{b}_j - Q^{b}_{j-1} \right)1_{ S^{b}_j = S^{b}_{j-1}} -   Q^{b}_{j-1 }  1_{ S^{b}_j < S^{b}_{j-1}} + Q^{b}_{j  }  1_{ S^{b}_j > S^{b}_{j-1}}, \nonumber\\
 V^{a}_j& = &\left( Q^{a}_j - Q^{a}_{j-1} \right)1_{ S^{a}_j = S^{a}_{j-1}} -   Q^{a}_{j-1 }  1_{ S^{a}_j > S^{a}_{j-1}} + Q^{a}_{j  }  1_{ S^{a}_j < S^{a}_{j-1}}.
 \end{eqnarray}
 We will describe the dynamics of the simplified order book using the process  $(X, Y)$ where
\begin{eqnarray}\label{main_processes}
X_j &=& V^{a}_j1_{  \textrm{$j$=ask-event} } - V^b_j1_{  \textrm{$j$=bid-event} },\notag\\
Y_j&=& Q^a_j 1_{  \textrm{$j$=ask-event} }- Q^b_j 1_{  \textrm{$j$=bid-event} }.
\end{eqnarray}

The modeling in this paper is based on the assumption in \eqref{a2} stated above, i.e. we assume that the best bid as well as the best ask price can only move by at most one tick (up or down) when an event occurs and  we first discuss two consequences of \eqref{a2}.  Firstly, assume the difference between the best bid  and the best ask prices is exactly one tick. Then these prices can only move if either the best bid queue $Q^b$ or the best ask queue $Q^a$ is depleted (i.e. one component of the order flow (\ref{order_flow}) hits 0) and when this occurs
the spread widens to two or more ticks. However, the assumption in \eqref{a2} now implies that at this event the spread is only allowed to widen to exactly two ticks. In particular, this means, for example, that if  best bid queue $Q^b$ is depleted then it is assumed that the second best bid queue,
one tick away, is non-empty. Secondly, assume the difference between the best bid price and the best ask price is at least two ticks. Then two things can happen: either the best bid queue $Q^b$ or the best ask queue $Q^a$ is depleted or a new limit order is placed between the best bid and ask. In the first case, based on \eqref{a2} and as above, the spread widens with exactly one tick. In the second case \eqref{a2} implies that the spread narrows by exactly one tick
since the new limit order is assumed to be submitted either one tick to the right of the best bid or one tick to the left of the best ask. In particular,
an other way of formulating \eqref{a2} is to state that
   \begin{eqnarray}\label{a4}\mbox{the bid/ask spread can only change by exactly one tick at each event.}
  \end{eqnarray}
However, we emphasis that we put no restrictions on the absolute size of the bid/ask spread.

Based on \eqref{a2}, \eqref{a4},  the dynamics of $(Q^b, Q^a, V^b, V^a)$ can be fully recovered from the dynamics of the process $(X, Y)$. Indeed, assume first the spread after event $j-1$ is exactly one tick. Then either the spread remains the same or it widens by one tick. If event $j$ occurs on the ask side of the order book then $Y_j= Q^a_j$ and if event $j$ occur on the bid side of the order book then $Y_j= -Q^b_j$. Hence $Y_j>0$ if and only if  event $j$ occurs on the ask side of the order book and  $Y_j<0$ if  event $j$ occurs on the bid side of the order book. Note that  the event $Y_j=0$ can not occur as, by definition, neither $Q^a$ nor $Q^b$ can hit zero. Furthermore, using this notation we see that $V^a_j = X_j 1_{Y_j >0}$ and $V^b_j = -X_j 1_{Y_j <0}$, while $Q^b$ and $Q^b$ are given in terms of $V^a$ and $V^b$ through the order flow equation (\ref{order_flow}) until the next price move. In this case, the next price move occurs
when the spread widens by one tick and then $Q^a$ ($Q^b$) is re-initialized by $Y$ ($-Y$) if the price moves up (down). Next, assume that
the spread after event $j-1$ is at least two ticks. Then either the spread remains the same or it widens by one tick or it narrows by one tick. In the first two cases we have as above that
$V^a_j = X_j 1_{Y_j >0}$ and $V^b_j = -X_j 1_{Y_j <0}$ while again $Q^b$ and $Q^b$ are given in terms of $V^a$ and $V^b$ through the order flow equation (\ref{order_flow}), until the next price move, and when this occurs $Q^a$ ($Q^b$) is re-initialized by $Y$ ($-Y$) as above. In the third case a new limit order is submitted either one tick to the right of the best bid or one tick to the left of the best ask. Assume that the new limit order is an ask order, i.e. an ask event occurs. Then $V_j^a>0$ and $X_j=V_j^a=Q_j^a=Y_j$ since
$V_j^a$ becomes the new best ask queue. Similarly, if the new limit order is a buy order, i.e. a buy event occurs, then $V_j^b>0$ and $X_j=-V_j^b=-Q_j^b=Y_j$ since
$V_j^b$ becomes the new best bid queue. We can conclude that $(V^b, V^a)$ can be fully recovered from the process $(X, Y)$ and hence  the dynamics of $(Q^b, Q^a)$ is determined by the dynamics of $(V^b, V^a)$ and hence by
$(X, Y)$. Finally, the dynamics of $(S^b, S^a)$ follows from $(X, Y)$ through $(Q^b, Q^a, V^b, V^a)$ and \eqref{a2}. In particular, the dynamics of $( S^b, S^a)$ can be described as
\[
[ S^b, S^a ] \longrightarrow  \left\{  \begin{array}{ll} ( S^b, S^a+ \delta ) \longrightarrow \left\{  \begin{array}{ll} ( S^b, S^a +2 \delta)     \\ ( S^b+\delta, S^a + \delta)\\ ( S^b-\delta, S^a + \delta)  \\ ( S^b, S^a )  \end{array} \right. \longrightarrow \dots    \\ ( S^b-\delta, S^a )  \longrightarrow \left\{  \begin{array}{ll} ( S^b-\delta, S^a + \delta)     \\ ( S^b , S^a  )\\ ( S^b-\delta, S^a - \delta)  \\ ( S^b-2 \delta, S^a )  \end{array} \right. \longrightarrow \dots \end{array}   \right.
\]
where $\delta$ is the tick. In conclusion,  based on \eqref{a2} the dynamics of $(S^b,S^a,Q^b,Q^a,V^b,V^a)$ is completely determined by $(X, Y)$. In particular, the complexity in our set up is reduced
to the process $(X, Y)$.

\subsection{Price dynamics} In the following we let 0 denote the index of the last event up to now and we let $\mathcal{F}_0$ be the history ($\sigma$-algebra) of all events that have occurred so far. We here discuss the dynamics of the best bid and ask prices and when we simply refer to the price we mean the mid price. The mid price can only change at an event as a consequence of that either the best bid and ask price change by the event.

Assume first that the bid-ask spread is exactly one tick, i.e. $S^a_0 = S^b_0 + \delta$. Define $l^{ b }_0 $ and $l^{a  }_0$ through
\begin{equation}\label{first_passage}
l^{b }_0 \equiv l^{b }_0 (Q^{b}_0) := \min \left\{ l \geq 1~ : ~ \sum_{j=1}^l V^{b}_{j} \le - Q^{b}_{0} \right\},
\end{equation}
and
\begin{equation}\label{first_passageA}
l^{a  }_0 \equiv l^{a  }_0 (Q^{a }_0) := \min \left\{ l \geq 1~ : ~ \sum_{j=1}^l V^{a }_{j} \le - Q^{a }_{0} \right\}.
\end{equation}
The next change in the mid price, best bid price or best ask price, occurs at the event index $\tau_0 = \min (l^{ b }_0,l^{a  }_0)$, i.e. when the best bid queue or the best ask queue is depleted, and if $l^{b }_0> l^{a }_0$, i.e. $\tau_0=l^{a }_0$,  then the next price move will be up and if $l^{b }_0< l^{a }_0$, i.e. $\tau_0=l^{b }_0$, then the next price move will be down. Naturally one can also ask what the probabilities for these events are. When $Q^{b,0}_0 \ll Q^{a,0}_0$, the bid queue will most likely be depleted before the ask queue, which means that the next price move will most likely be down. However, the outcome is less obvious when $Q^{b}_0 $ is close to $  Q^{a}_0$. Let
\begin{eqnarray}
p^+_0(q^b,q^a) &:=& \mathbb{P} \left( l^{a }_0 (q^a) < l^{b }_0 (q^b) | \; \mathcal{F}_0 \right),\notag \\
 p^-_0(q^b,q^a)&:=& \mathbb{P} \left( l^{a }_0 (q^a) > l^{b }_0 (q^b)  | \; \mathcal{F}_0\right),
\end{eqnarray}
whenever $(q^b,q^a)\in\mathbb R^2$.  Then $p^+_0(q^b,q^a)$  is the probability that ask queue, currently having size $q^a$, will be depleted
before that  bid queue, currently having size $q^b$. In particular, conditioned on $(q^b,q^a)$, $p^+_0(q^b,q^a)$ is the probability that the next prove move will
be up. $p^-_0(q^b,q^a)$ have a similar interpretation. By (\ref{first_passage}), we see that $ q^a\longmapsto l^{a }_0 (q^a)$ and $ q^b\longmapsto l^{b }_0 (q^b)$ are non decreasing functions. Furthermore,  for $q^b>0$ fixed  the function $ q^a \longmapsto p^+_0(q^b,q^a)$ is  non increasing and the function $ q^a \longmapsto p^-_0(q^b,q^a)$ is  non decreasing. Also
\begin{eqnarray}
&&\lim_{q^a\to 0^+}p^+_0(q^b,q^a)=1,\ \lim_{q^a\to 0^+}p^-_0(q^b,q^a)=0,\notag\\
&&\lim_{q^a\to \infty}p^+_0(q^b,q^a)=0,\ \lim_{q^a\to \infty}p^-_0(q^b,q^a)=1.
\end{eqnarray}
In particular, given $q^b>0$, there exists a unique point $q^\ast \equiv q^\ast(q^b)$ such that $  p^+_0(q^b,q^a)\leq p^-_0(q^b,q^a) $ for every $q^a\geq  q^\ast(q^b)$ and $ p^+_0(q^b,q^a) >  p^-_0(q^b,q^a) $ for every $q^a < q^\ast(q^b)$. Note that $q^\ast(q^b)$ can be written as
\begin{equation}
q^\ast(q^b) = \inf \left\{ q^a\in \mathbb{R} ~ : ~   p^+_0(q^b,q^a) \geq  p^-_0(q^b,q^a)\right\} \vee 0.
\end{equation}
Furthermore, the function $q^b \longmapsto q^\ast(q^b)$ is non decreasing. The model predicts a price increase if $Q^{a}_0< q^\ast (Q^{b}_0) $ and a price decrease if $Q^{a}_0> q^\ast (Q^{b}_0)$.

 Next, using $\tau_0$ as introduced above, we see that if $S^{a}_0 = S^{b}_0 + \delta$, i.e. the spread is exactly one tick, then the state of $(S^{a(b)}, Q^{a(b)})$ after event $  k  \in \left\{ 0, \; \dots,   \tau_0 -1 \right\}$ is given by
\begin{equation}
 S^{a}_k= S^{a}_0 ,~ S^{b}_k = S^{b}_0, ~ Q^{a}_k = Q^{a}_0 + \sum_{j= 1}^k V^{a}_{j},~ Q^{b}_k = Q^{b}_0 + \sum_{j= 1}^k V^{b}_j.
  \end{equation}
  When either the best bid or ask queue is depleted, i.e at $\tau_0  $, the spread widens immediately by one tick and $S_0^{a(b)} , ~ Q_0^{a(b)}$ are re-initialized according to
  \begin{eqnarray}\label{spread_up}
 S^{a}_{\tau_0} &=& S^{a}_0 + \delta  \; 1_{l^{a }_0 <l^{b }_0}, ~ S^{b}_{\tau_0} = S^b_0 - \delta  \; 1_{l^{a }_0 >l^{b }_0}, \nonumber\\
 Q^{a}_{\tau_0} &=& \left( Q^{a}_{ j}  + \sum_{j= 1}^{\tau_0} V^{a}_{j} \right)1_{l^{a }_0 >l^{b }_0} + Y_{ \tau_0  } \; 1_{l^{a }_0 <l^{b }_0}, \nonumber\\
 Q^{b}_{\tau_0} &=& \left( Q^{b}_{ 0}  + \sum_{j= 1}^{\tau_0} V^{b}_{j} \right)1_{l^{a }_0 <l^{b }_0}  -  Y_{ \tau_0 } \;1_{l^{a }_0 >l^{b }_0}.
\end{eqnarray}

By the above the spread is now two ticks and again are now two scenarios for the next change in the spread. Either it widens again by one tick, in analogy with the outline above, or it narrows again if a new limit order in placed between the best bid and ask. Let
\begin{equation}\label{l_a_b}
l^{a(b) }_1  := \min \left\{ l \geq \tau_0 + 1~ : ~ \sum_{j=\tau_0 + 1}^l V^{a(b)}_{j} \le - Q^{a(b)}_{\tau_0} \right\}
\end{equation}
and  \begin{eqnarray}\label{tilde_l_a_b}
\tilde{l}^{a}_1 &=&  \min \left\{ l \geq \tau_0 +1~ :\mbox{ a new limit ask order is submitted between}\right\},\notag\\
\tilde{l}^{b}_1 &=&  \min \left\{ l \geq \tau_0 +1~ :\mbox{ a new limit bid order is submitted between}\right\}.
\end{eqnarray}
The time $l=\tilde{l}_1^{a } \wedge \tilde{l}_1^{b }$ is also the first event time, after $\tau_0$, after which  $X_l$ and $Y_l$ have the same sign and $|Y_l| \le |X_l|$. Indeed,  assume first the price does not change after event $l$, then $Y_l =  - Q^b_{l-1} + X_l  $Ê if $ Y_l< 0 $ and $Y_l =   Q^a_{l-1} +X_l  $ Êif $ Y_l>0 $. Hence either $Y_l$ and $X_l$ do not have the same sign or $|X_l| < |Y_l|$. Next assume the spread widens after event $l$, then either $Y_l > 0$, in which case the best ask queue is depleted and $ X_l = V^a_l = - Q^a_{l-1} <0$, or  $Y_l < 0$ and $ X_l = -V^b_l =   Q^b_{l-1} >0$. Hence $X_l$ and $Y_l$ do not have the same sign.  Finally, assume the spread narrows after event $l$, then a new limit buy (or sell) order is submitted between the best bid price and the best ask price and $X_l = V^a_l = Q^a_l = Y_l $, if $Y_l>0$ (new sell order), and $X_l = -V^b_l = -Q^b_l = Y_l,$ if $Y_l<0$ (new buy order). We can therefore write $\tilde{l}_1^{a }$ and $\tilde{l}_1^{b}$  as
\begin{eqnarray}
\tilde{l}_1^{b} &=& \min\left\{ l \geq \tau_0 +1 ~: ~ X_l  \le  Y_l  <0  \right\}, \nonumber\\
\tilde{l}_1^{a }& =& \min\left\{ l \geq \tau_0 +1 ~: ~  0 < Y_l \le X_l \right\}.
\end{eqnarray}
Finally, let \begin{equation}
\tau_1 = l^{a  }_1 \wedge l^{b }_1 \wedge \tilde{l}^{b}_1 \wedge \tilde{l}^{a}_1.
\end{equation}
Using this notation, we see that the next price move, after $ \tau_0$, occurs at event $\tau_1$. In particular, if spread widens again by one tick
then $\tau_1 = l^{a  }_1 \wedge l^{b }_1 $ and we see that we can describe the dynamics of $(S^b,S^a,Q^b,Q^a)$ as above. Similarly, if the spread next narrows by one tick then $\tau_1 =  \tilde{l}^{b}_1 \wedge \tilde{l}^{a}_1$ and once the order is submitted the size of the new queue  is either $V^a_{\tilde{l}^{a}_1}$ or $V^b_{\tilde{l}^{b}_1}$. Under this scenario, at $ \tau_1$, the best bid and ask prices move from $ ( S^{b}_{\tau_0 }, S^{a}_{\tau_0} ) $ to $ ( S^{b}_{\tau_1 }, S^{a}_{\tau_1} ) $, where
\begin{eqnarray}\label{queue_1}
 S^{b}_{\tau_1} &=&   S^{b}_{\tau_0} + \delta (  1_{ \tau_1 = \tilde{l}^b_1} -  1_{ \tau_1 = l^b_1 } ), \nonumber\\
 S^{a}_{\tau_1} &=&  S^{a}_{\tau_0} + \delta (  1_{ \tau_1 = l^a_1 } - 1_{ \tau_1 = \tilde{l}^a_1}   ),
\end{eqnarray}
and the queues $ Q^{a}_{\tau_1}$ and $ Q^{b}_{\tau_1}$ are re-initialized to
 \begin{eqnarray}\label{queue_1_2}
 Q^{b}_{\tau_1} &=&  V^{b }_{\tau_1} \; 1_{ \tau_1 =\tilde{l}^{b}_1  } -Y_{\tau_1} 1_{ \tau_1 =l^b_1   }+
 \left( Q^{b}_{\tau_0}+  \sum_{l=\tau_0 + 1}^{\tau_1 } V^{b}_l \right)\; 1_{ \tau_1  < \tilde{l}^{b}_1 \wedge l^{b}_1  }, \nonumber\\
 Q^{a}_{\tau_1} &=&  V^{a }_{\tau_1} \; 1_{ \tau_1 =\tilde{l}^{a}_1  } + Y_{\tau_1} 1_{ \tau_1 =l^a_1   }+
 \left( Q^{a}_{\tau_0}+  \sum_{l=\tau_0 + 1}^{\tau_1 } V^{a}_l \right)\; 1_{ \tau_1  < \tilde{l}^{a}_1 \wedge l^{a}_1  }.
\end{eqnarray}

We emphasize that while we have here used a mix of the notation  $(S^b,S^a,Q^b,Q^a,V^b,V^a)$ and  $(X, Y)$, each component
of $(S^b,S^a,Q^b,Q^a,V^b,V^a)$ can in the above displays be expressed solely using $(X, Y)$ as discussed in the previous subsection.

\end{document}